\documentclass[oneside,german,english]{amsart}
\usepackage{palatino}
\usepackage[T1]{fontenc}
\usepackage[latin1]{inputenc}
\usepackage{amssymb}
\IfFileExists{url.sty}{\usepackage{url}}
                      {\newcommand{\url}{\texttt}}

\makeatletter

 \theoremstyle{plain}
 \theoremstyle{plain}    
 \newtheorem*{thm*}{Theorem} 
 \theoremstyle{plain}    
 \newtheorem*{lem*}{Lemma} 

\ifx\hyperlink\undefined
\newcommand{\hrlb}{ }

\newcommand{\hyperlink}[2]{#2}
\else
\newcommand{\hrlb}{\\}
\fi
\def\UrlBigBreaks{\do\.\do\?\do\-\do\:\do@url@hyp}

\usepackage{babel}
\makeatother
\begin{document}

\title{Percolation, Perimetry, Planarity}

\author{Gady Kozma}

\address{Institute for Advanced Study, 1 Einstein drive, Princeton, NJ 08540,
USA.}

\email{gady@ias.edu}

\begin{abstract}
Let $G$ be a planar graph with polynomial growth and isoperimetric
dimension bigger than $1$. Then the critical $p$ for Bernoulli percolation
on $G$ satisfies $p_{c}<1$.
\end{abstract}
\maketitle
This short note relates to a famous problem first posed in the seminal
paper \cite{BS96}: find some general conditions on a graph $G$ under
which $p_{c}<1$. In particular they conjectured (question 2 ibid.)
that if the isoperimetric dimension is $>1$ then this holds. The
isoperimetric dimension is defined for an infinite connected graph
$G$ by \[
\dim(G):=\sup\left\{ d\geq1:\inf_{S}\frac{|\partial S|}{|S|^{(d-1)/d}}>0\right\} \]
where the infimum is over all finite non-empty sets of vertices $S$.
$\partial S$ is the set of edges with one vertex in $S$ and the
other in $G\setminus S$.

A number of partial cases of this conjecture have been proved. In
\cite[theorem 2]{BS96} it was proved that a positive Cheeger constant
(i.e.~$|\partial S|>c|S|$) implies $p_{c}<1$. In \cite{PS04} this
was proved for graphs which satisfy some complicated conditions on
the geometry of the minimal cut-sets. Here we shall show this under
one technical condition (polynomial growth) and one more essential
condition (planarity).

I wish to thank Itai Benjamini for many fruitful discussions. G\'abor
Pete, Paul Seymour and Carsten Thomassen helped with references and
advice.

\begin{thm*}
Let $G$ be a planar graph with no vertex accumulation points such
that
\begin{enumerate}
\item There exist numbers $K$ and $D$ such that for all $v\in G$ and
any $r\geq1$ one has for the open ball $B(v,r)$ that the number
of vertices satisfies $|B(v,r)|\leq Kr^{D}$; and
\item There exist numbers $k,\epsilon>0$ such that for any finite non-empty
set of vertices $S$, $|\partial S|\geq k|S|^{\epsilon}$.
\end{enumerate}
Let $p_{c}$ be the critical $p$ for independent bond percolation
on $G$. Then $p_{c}<1$.
\end{thm*}
A vertex accumulation point is a point $x\in\mathbb{R}^{2}$ such
that every neighborhood of $x$ contains infinitely many vertices
of $G$. To be more precise, the theorem holds whenever the abstract
planar graph $G$ has a plane representation $\Gamma$ with no vertex
accumulation points. We assume $G$ has no loops or multiple edges.

\begin{proof}
Let us first dispose of uninteresting cases resulting from stating
the theorem in too much generality. We may assume every vertex of
$G$ has only finitely many neighbors since otherwise $p_{c}=0$.
Further, By Wagner's theorem \cite{W67} a planar graph contains no
more than a countable number of vertices with degree $\geq3$. On
the other hand, our graph $G$ cannot contain a component with all
degrees $\leq2$ since that would be a finite or infinite line and
would violate the assumption that the isoperimetric dimension is $>1$.
With Wagner's theorem we get that $G$ is countable. Hence the existence
of an infinite cluster is a measurable event, $p_{c}$ is well defined
and our theorem makes sense. This allows us to assume $G$ is connected
by restricting to any infinite component.

We use Peierls argument, which states that $p_{c}<1$ whenever the
number of minimal cut-sets of some vertex $v\in G$ of size $n$ is
exponential (or less) in $n$. See e.g.~\cite[page 16]{G99}. Here
a cut-set of $v$ is a set of edges $S$ such that $v$ is in a finite
component of $G\setminus S$. The cut-set is minimal if every $S'\subsetneqq S$
is not a cut-set. Since $G$ is planar we can define its dual graph
$G^{*}$ by making any face of $G$ to be a vertex of $G^{*}$ and
matching any edge $e$ of $G$ to an edge $e^{*}$ of $G^{*}$ between
the two faces of $G$ on the two sides of $e$. Note that $G^{*}$
may contain loops and multiple edges. Formally we require the following
from $G^{*}$:
\begin{enumerate}
\item There is a one-to-one onto correspondence between $E(G)$, the edge
set of $G$ and $E(G^{*})$, which we denote by $*$.
\item A minimal cut-set of $G$ is carried by $*$ to a cycle (i.e.~a simple
closed path) of $G^{*}$ and vice versa.
\end{enumerate}
The existence of $G^{*}$ is well known for finite graphs but for
infinite graphs I was not able to find a completely suitable reference,
hence we shall show the existence of $G^{*}$ in the appendix.

Thus we need to show that the number of cycles of length $n$ separating
$v^{*}$ from infinity is $\leq C^{n}$ for some $n$. If, for example,
$G^{*}$ happens to be a graph with bounded degree (as in \cite{PS04})
we could have finished here. In general, however, $G^{*}$ could have
unbounded degree and might not even be locally finite (a graph is
locally finite if every vertex has finite degree).

It will be easier to examine simple paths. Let us therefore fix two
vertices $a^{*}\neq b^{*}\in G^{*}$ and examine the number $p=p(a^{*},b^{*};n)$
of simple paths of length $n$ starting from $a$ and ending at $b$.
We shall now show that $p(a^{*},b^{*};n)\leq C^{n}$ and then we shall
be mostly finished. Let \[
p(n):=\max_{a^{*}\neq b^{*}\in G^{*}}p(a^{*},b^{*};n).\]
 We wish to find an inequality connecting $p(n)$ and $p(2n)$.

Let therefore $a^{*}\neq b^{*}\in G^{*}$ and let $\gamma^{*}$ be
one simple path of length $2n$ between $a^{*}$ and $b^{*}$ which
we will use as a reference. Let $\delta^{*}$ be a second such path
--- we wish to bound the number of possibilities for $\delta^{*}$.
For every edge $e^{*}$ of $\gamma^{*}$ there are two vertices of
$G$ adjacent to $e$ --- let $\{ v_{1},\dotsc,v_{4n}\}$ be the complete
list. Let $f^{*}$ be the $n$'th edge of $\delta^{*}$, and let $w_{1}$,
$w_{2}$ be the two vertices of $G$ on both sides of $f$. Since
$\gamma^{*}\cup\delta^{*}$ is a closed path, either $f^{*}\in\gamma^{*}$
or there exists a cycle $\beta^{*}\subset\gamma^{*}\cup\delta^{*}$
containing $f^{*}$. Assume the second. Since $\beta$ is a minimal
cut-set in $G$ we have that it is the boundary of some finite connected
set $Q$ and therefore one of the $w_{j}$-s is in $Q$. Denote it
by $w$. It is not possible for $\beta^{*}\subset\delta^{*}$ (since
$\delta^{*}$ is simple) so there exists at least one of the $\{ v_{1},\dotsc,v_{4n}\}$
which is in $Q$. Denote it by $v$. Since $|\partial Q|=|\beta|\leq|\gamma^{*}\cup\delta^{*}|\leq4n$
then by the isoperimetric inequality we see that $|Q|\leq(4n/k)^{1/\epsilon}$.
Since $Q$ is connected the distance in $G$ between $v$ and $w$
must be $\leq(4n/k)^{1/\epsilon}$. We think about this as $w\in\bigcup_{i}B(v_{i},(4n/k)^{1/\epsilon})$
and we have also subsumed the case that $f^{*}\in\gamma^{*}$ which
wasn't covered by the argument above (choose $w$ from $w_{1}$, $w_{2}$
arbitrarily in this case). Polynomial growth shows that $w$ has $\leq4nK(4n/k)^{D/\epsilon}$
possibilities. Using the polynomial growth inequality with $r=2$
we know that every vertex of $G$ has no more than $K2^{D}$ adjacent
edges, so $f$ has $\leq4nK^{2}(8n/k)^{D/\epsilon}$ possibilities
and the same holds (with another factor of $2$) for the middle vertex
of $\delta^{*}$. Since this holds for every $\delta^{*}$ we are
left with the inequality\[
p(2n)\leq8nK^{2}(8n/k)^{D/\epsilon}p(n)^{2}.\]
An identical argument shows that $p(1)\leq(2/k)^{1/\epsilon}$ and
from here a simple induction shows that $p(2^{n})\leq\exp(C2^{n})$
where $C=C(k,K,\linebreak[0]\epsilon,D)$ is independent of $n$.
This obviously implies $p(n)\leq C^{n}$ for some other $C$.

To finish the theorem, let $v\in G$ and examine a minimal cut-set
$\gamma$ of size $n$. Since $\gamma$ is the boundary of the component
of $G\setminus\gamma$ containing $v$, we see that this component
has size $\leq(n/k)^{1/\epsilon}$ and in particular it is contained
in $B(v,(n/k)^{1/\epsilon})$. By polynomial growth we see that there
are no more than $K(n/k)^{D/\epsilon}$ possible vertices in the component
and hence no more than $K^{2}(2n/k)^{D/\epsilon}$ edges which may
participate in the cut-set. For each of these edges $e$, a cut-set
containing $e$ is, viewed in $G^{*}$, $e^{*}\cup\{$an open simple
path of length $n-1\}$. Therefore the number of possibilities for
a minimal cut-set is \[
\leq K^{2}(2n/k)^{D/\epsilon}p(n-1)\leq C^{n}\]
for some $C$. This proves the theorem.
\end{proof}

\subsection*{Appendix: duality}

Duality of infinite planar graphs was investigated by Thomas\-sen
\cite{T80,T82} and recently by Bruhn and Diestel \cite{BD06} but
we cannot use either for the following reason. Both define the dual
(denote it by $\widehat{G}$ to differentiate) such that a cycle in
$\widehat{G}$ corresponds to a bond in $G$ where a bond is a minimal
set dividing $G$ into components \emph{but without requiring any
one of them to be finite}. Take as an example $G$ to be an infinite
bi-directional line. $G^{*}$ (which is, in this case, uniquely defined
up to the addition of isolated vertices) is a graph with two vertices
connected by infinitely many edges, while $\widehat{G}$ is a single
vertex with infinitely many loops. Hence our minimal cut-sets correspond
in $\widehat{G}$ to a not-necessarily simple closed path and the
proof does not go through. 

It is quite likely that the condition that the graph has no accumulation
points is redundant (with the same proof structure). Condition (2)
in the definition of $G^{*}$ must be relaxed by removing the {}``vice
versa'' claim, but this can be worked around in the proof of the
theorem. Hence the only obstacle is the existence of $G^{*}$. It
would be interesting to show that $G^{*}$ exists under, say, the
condition $(*)$ of \cite{BD06}.

\begin{lem*}
Let $G$ be a connected locally finite planar graph with no accumulation
points. Then a $G^{*}$ satisfying the requirements (1) and (2) exists.
\end{lem*}
\begin{proof}
[Proof sketch]By \cite[theorem 3]{T77} there exists a straight line
triangulation $\Delta$ such that $G$ is a subgraph of $\Delta$.
Define two triangles of $\Delta$ to neighbor if they have a common
edge and it is not in $G$. We will call the components of the neighborhood
graph on the triangles of $\Delta$ {}``faces of $G$'' and they
will serve as the vertex set for $G^{*}$. For every edge $e$ of
$G$ we will define the edge $e^{*}$ to be adjacent to the faces
containing its two neighboring triangles (which may be the same).
Hence we need only show the relation between minimal cut-sets of $G$
and cycles of $G^{*}$.

To show this we will construct $G^{*}$ as a geometric dual of $G$.
In any face $f$ of $G$ pick an arbitrary edge $e$ in $\partial f$,
let $T$ be the triangle of $f$ containing $e$ and let $p\in\mathbb{R}^{2}$
be the middle $T$. Now trace around $\partial f$ starting from $T$
and $e$ using the triangulation structure of $\Delta$ i.e.~at any
vertex turn around until hitting an edge of $G$ and continue to the
appropriate triangle. This process reaches the entire boundary of
$f$. To see this, use the tracing process to construct a simple path
$\rho\subset\mathbb{R}^{2}$ ($\partial f$ need not be simple) slightly
inside $f$ (e.g.~in every triangle $S$ of $\Delta$ participating
in the tracing process make $\rho$ go through the third of $S$ closest
to $G$). If the resulting $\rho$ is infinite it must go to infinity
(because $G$ has no accumulation points) and we can add the point
at infinity to close it. Let now $e'\in\partial f$ be arbitrary and
we wish to show that it participated in the tracing process. Let $S_{1},\dotsc,S_{n}$
be a path of neighboring triangles such that $e'\in\partial S_{1}$
and $S_{n}=T$. Let $\sigma\subset\mathbb{R}^{2}$ be a path starting
in the middle of $e'$, linear in every $S_{i}$ and passing through
the middle of the edge joining $S_{i}$ and $S_{i+1}$, and finally
ending in $p$. We use Jordan's theorem with $\rho$ and get two components
$A$ and $B$ and examining $T$ it is easy to see that $p$ is in
one (say $A$) while $e$ (and therefore, by connectivity, of all
$G$) is in the other. Therefore $\rho\cap\sigma\neq\emptyset$. But
by definition this can only happen in $S_{1}$ and therefore $\rho$
passed through $S_{1}$ and therefore $e'$ participated in the tracing
process.

The tracing process allows to connect $p$ to the middle of any $e'\in\partial f$
by disjoint curves. These curves will be {}``half edges'', the other
half coming from the other face. Hence we have completed a description
of $G^{*}$ as a planar graph. Further, $e\cap f^{*}\neq\emptyset$
iff $e=f$ and then they intersect at a single point and are transversal
at that point. This immediately shows that any cycle $\gamma^{*}$
in $G^{*}$ maps to a cut in $G$ --- Jordan's theorem shows that
the set $Q$ of vertices of $G$ inside $\gamma^{*}$ is disjoint
from the ones outside $\gamma^{*}$, and the transversality shows
that for every $e^{*}\in\gamma^{*}$ exactly one of the vertices adjacent
to $e$ is in $Q$, so $Q\neq\emptyset$. The fact that there are
no accumulation points in $G$ shows that $Q$ is finite and hence
$\gamma$ is a cut. 

On the other hand, if $\gamma$ is a minimal cut then every vertex
in $\gamma^{*}$ must have degree at least $2$. To see this take
one $e^{*}\in\gamma^{*}$ and let $p^{*}$ be an adjacent face of
$G$, and assume no other $f^{*}\in\gamma^{*}$ is adjacent to $p^{*}$.
Let $\delta$ be the boundary of $p^{*}$ which is either a closed
path or a bi-directionally infinite path in $G$ (in neither case
necessarily simple). If $\delta$ is finite then it forms a path in
$G\setminus\gamma$ between the two ends of $e$, in contradiction
of the minimality of $\gamma$. If $\delta$ is infinite then both
ends of $e$ are connected to an infinite number of vertices in $G$,
in contradiction to the fact that one side must be finite.

Now, a graph $\gamma^{*}$ with minimal degree $\geq2$ must contain
a cycle, $\delta^{*}$. By the previous argument, $\delta$ would
be a cut. But since we assumed $\gamma$ is minimal, we must have
$\gamma=\delta$ and hence $\gamma^{*}$ is a cycle. Conversely, if
$\gamma^{*}$ is a cycle then $\gamma$ is a cut. If $\delta\subset\gamma$
is minimal then $\delta^{*}\subset\gamma^{*}$ is a cycle. But cycles
don't contain subcycles so $\delta^{*}=\gamma^{*}$ and we have established
both directions of the correspondence between minimal cuts of $G$
and cycles of $G^{*}$, and we are done.
\end{proof}


\begin{thebibliography}{BS96}
\bibitem[BS96]{BS96}Itai Benjamini and Oded Schramm, \emph{Percolation beyond $\mathbb{Z}^{d}$,
many questions and a few answers}, Electron. Comm. Probab. \textbf{1:8}
(1996), 71--82.\hrlb  \url{http://www.math.washington.edu/~ejpecp/ECP/viewarticle.php?id=1561&layout=abstract}
\bibitem[BD06]{BD06}Henning Bruhn and Reinhard Diestel, \emph{Duality in infinite graphs},
Combin. Probab. Comput. \textbf{15:1--2} (2006), 75--90. \hrlb \url{http://www.math.uni-hamburg.de/home/diestel/papers/Duality.pdf}
\bibitem[G99]{G99}Geoffrey Grimmett, \emph{Percolation}, second edition. \foreignlanguage{german}{Grundlehren
der Mathematischen Wissenschaften} {[}Fundamental Principles of Mathematical
Sciences{]}, 321. Springer-Verlag, Berlin, 1999.
\bibitem[PS04]{PS04}Aldo Procacci and Benedetto Scoppola, \emph{Infinite graphs with a
nontrivial bond percolation threshold: some sufficient conditions},
J. Statist. Phys. \textbf{115:3--4} (2004), 1113--1127.\hrlb  \url{http://www.mat.ufmg.br/~aldo/papers/ps7.pdf}
\bibitem[T77]{T77}Carsten Thomassen, \emph{Straight line representations of infinite
planar graphs}, J. London Math. Soc. (2) \textbf{16:3} (1977), 411--423.
\bibitem[T80]{T80}Carsten Thomassen, \emph{Planarity and duality of finite and infinite
graphs}, J. Combin. Theory Ser. B \textbf{29:2} (1980), 244--271.
\bibitem[T82]{T82}Carsten Thomassen, \emph{Duality of infinite graphs}, J. Combin. Theory
Ser. B \textbf{33:2} (1982), 137--160.
\bibitem[W67]{W67}Klaus Wagner, \emph{Fastpl\"attbare Graphen} (German), J. Combinatorial
Theory \textbf{3} (1967) 326--365.
\end{thebibliography}
\end{document}